\documentclass[12pt]{amsart}

\usepackage{amsthm,amssymb,amsfonts,amsmath,enumerate,graphicx,color,tikz}
\usepackage{mathrsfs}
\usepackage[margin=1in]{geometry}

\usepackage{caption}
\usepackage{subcaption}

\usepackage[utf8]{inputenc}

\usepackage{hyperref}

\newcommand{\Pc}{\mathcal{P}}

\newcommand{\Cc}{\mathcal{C}}
\newcommand{\Hc}{\mathcal{H}}

\newcommand{\Pb}{\mathbf{P}}
\newcommand{\Ib}{\mathbf{I}}

\newcommand{\C}{\mathbb{C}}
\newcommand{\Z}{\mathbb{Z}}
\newcommand{\N}{\mathbb{N}}
\newcommand{\R}{\mathbb{R}}
\newcommand{\Q}{\mathbb{Q}}

\newcommand{\des}{\mathrm{des}}

\newcommand{\cone}{\mathrm{cone}}
\newcommand{\conv}{\mathrm{conv}}

\renewcommand{\phi}{\varphi}
\renewcommand{\emptyset}{\varnothing}

\def\c{{\boldsymbol c}}
\def\d{{\boldsymbol d}}
\def\e{{\boldsymbol e}}

\def\p{{\boldsymbol p}}
\def\r{{\boldsymbol r}}
\def\s{{\boldsymbol s}}
\def\t{{\boldsymbol t}}
\def\u{{\boldsymbol u}}

\def\x{{\boldsymbol x}}

\newcommand\commentout[1]{}

\newcommand{\Znn}{\Z_{\geq 0}}
\newcommand{\Zpos}{\Z_{\geq 1}}

\newcommand{\Pbn}{\Pb_n^{(\s)}}
\newcommand{\Ibn}{\Ib_n^{(\s)}}

\newcommand{\asc}{\operatorname{asc}}
\newcommand{\Asc}{\operatorname{Asc}}

\newcommand{\Sn}{\mathfrak{S}_n}
\newcommand{\Lc}{\mathcal{L}}

\newtheorem{theorem}{Theorem}[section]
\newtheorem{corollary}[theorem]{Corollary}

\newtheorem{lemma}[theorem]{Lemma}
\newtheorem{conjecture}[theorem]{Conjecture}

\theoremstyle{remark}

\newtheorem{remark}[theorem]{Remark}
\theoremstyle{definition}
\newtheorem{definition}[theorem]{Definition}
\newtheorem{question}[theorem]{Question}
\newtheorem{problem}[theorem]{Problem}

\begin{document}

\title{Polyhedral geometry for  lecture hall partitions}

\author{McCabe Olsen}
\address{Department of Mathematics\\
         The Ohio State University\\
         231 W. 18th Ave., Columbus, OH 43210, USA}
\email{olsen.149@osu.edu}




\date{\today}

\thanks{This manuscript was
prepared as a contribution to the conference proceedings of the 2018 \emph{Workshop
on Lattice Polytopes} at Osaka University. The author thanks the organizers,
Takayuki Hibi and  Akiyoshi Tsuchiya, for their support during this workshop.
The author also thanks Benjamin Braun and Liam Solus for helpful comments and feedback.}


\begin{abstract}
Lecture hall partitions are a fundamental combinatorial structure which have been studied extensively over the past two decades. 
These objects have produced new results, as well as reinterpretations and generalizations of classicial results, which are of interest in combinatorial number theory, enumerative combinatorics, and convex geometry.
In a recent survey of Savage \cite{Savage-LHP-Survey}, a wide variety of these results are nicely presented.
However, since the publication of this survey, there have been many new developments related to the polyhedral geometry and Ehrhart theory arising from lecture hall partitions.
Subsequently, in this survey article, we focus exclusively on the polyhedral geometric results in the theory of lecture hall partitions in an effort to showcase these new developments.
In particular, we highlight results on lecture hall cones, lecture hall simplices, and lecture hall order polytopes.
We conclude with an extensive list of open problems and conjectures in this area. 
\end{abstract}

\maketitle



\section{Introduction}

For $n\in\Zpos$, let $[n]:=\{1,2,\ldots,n\}$.
Let $\s=(s_1,\ldots,s_n)\in\Zpos^n$ a sequence of positive integers.
 The set of \emph{$\s$-lecture hall partitions} is the set 
	\[
	L_n^{(\s)}:=\left\{ \lambda\in\Z^n \, : \, 0\leq \frac{\lambda_1}{s_1}\leq \cdots\leq \frac{\lambda_n}{s_n} \right\}.
	\]
In the case that $\s$ is weakly increasing sequence, $L_n^{(\s)}$ refines the set of partitions with at most $n$ parts.	
In other cases, one may actually have a composition rather than a partition, though much of the literature makes no distinction here.
The study of lecture hall partitions was initiated by Bousquet-M\'{e}lou and Eriksson in a pair of seminal papers \cite{BME-LHP1,BME-LHP2}.
In the two decades to follow, these objects have been studied in a variety of algebraic, combinatorial, geometric, and number theoretic mediums.
Much of this work is summarized in the excellent survey of Savage \cite{Savage-LHP-Survey}.

The goal of this survey is to provide an update to the survey of Savage \cite{Savage-LHP-Survey}. 
In particular, we focus exclusively on results involving convex geometric objects arising from lecture hall partitions taking into account a significant number of results appearing after the publication of the Savage survey.
The structure of this note is as follows. 
In Section \ref{sec:background}, we provide a brief review of the necessary polyhedral geometry, Ehrhart theory, and local Ehrhart theory, as well as define the convex geometric obejects arising from lecture hall partitions, namely \emph{lecture hall cones}, \emph{lecture hall simplices}, \emph{rational lecture hall simplices}, and \emph{lecture hall order polytopes}.
In Section \ref{sec:Cones}, we outline known results of lecture hall cones.
In Section \ref{sec:Simplices}, we discuss the known results regarding lecture hall simplices.
Section \ref{sec:order} focuses on results of lecture hall order polytopes. 
Finally, Section \ref{sec:Open} outlines some of the remaining open problems in the polyhedral geometry of lecture hall partitions.

\section{Background}\label{sec:background}

In this section, we provide the necessary background for the results on the polyhedral geometry of lecture hall partitions. 
First, we review the standard tools and definitions in Ehrhart theory and polyhedral geometry.
We then discuss the notions of local Ehrhart theory and box unimodal triangulations.
Finally, we briefly introduce the polyhedral geometric objects arising for lecture hall partitions. 
Subsequently, some or all of these sections may be safely skipped by the experts.

\subsection{Classical Ehrhart Theory and Polyhedral Geometry}

We recall a few definitions from polyhedral geometry.
A \emph{polyhedral cone} $C$ in $\R^n$ is the solution set to some finite collection of linear inequalities $Ax\geq 0$ for some real matrix $A$, or equivalently given some elements $w_1,w_2,\cdots,w_j \in \R^n$, 
	\[
	C=\operatorname{span}_{\R_{\geq 0}}\{w_1,w_2,\cdots,w_j  \}.
	\]
The elements $w_i$ are called \emph{ray generators}. The cone $C$ is said to be \emph{rational}	if the matrix $A$ contains rational entries (equivalently if each $w_i\in\Q^n$), it is said to be \emph{simplicial} if it is defined by $n$ independent inequalities (equivalently if $j=n$ and $\{w_i\}_{i=1}^n$ are linearly independent),  and it is said to be \emph{pointed} if it does not contain a linear subspace of $\R^n$. Let $C^\circ$ denote the interior of $C$. 

If $C\subset\R^n$ is a pointed, rational cone, a \emph{proper grading} of $C$ is a function $g:C\cap\Z^n\to \Znn^r$, for some $r$, satisfying (i) $g(\lambda+\mu)=g(\lambda)+g(\mu)$; (ii) $g(\lambda)=0$ implies $\lambda=0$; and (iii) for any $v\in \N^r$, $g^{-1}(v)$ is finite. Moreover, the integer points $C\cap \Z^n$ form a semigroup. 
Semigroups of this type have unique minimal generating sets known as the \emph{Hilbert basis} of $C$. 
Additionally, the lattice points in a pointed rational cone give rise to a semigroup algebra structure. We call this the \emph{affine, normal semigroup algebra of $C$}, which is denoted $\C[C]:=\C[C\cap\Z^n]$
For background and details see \cite{BeckRobins-CCD,MillerSturmfels-CCA}.

We say that a pointed, rational cone $C\subset\R^n$ is  \emph{Gorenstein} if there exists a point $\c\in C^\circ$ such that $C^\circ \cap \Z^n=\c+(C\cap \Z^n)$. This point is known as the \emph{Gorenstein point} of $C$. Due to theorems of Stanley \cite{Stanley-HilbertFunctions}, this notion of Gorenstein is equivalent to the commutative algebra notion of Gorenstein, as $C$ is Gorenstein if and only if the algebra $\C[C]$ is Gorenstein. For reference and commutative algebra details, see \cite{BrunsHerzog,StanleyGreenBook}.

We say that $\Pc\subset \R^n$ is a \emph{convex polytope} of dimension $d\leq n$ if it can be expressed as the convex hull of finitely many points in  
	\[
	\Pc=\conv\{v_0,\ldots,v_m \, : \, v_i\in \R^n\}:=\left\{\sum_{i=0}^m \mu_i v_i \, : v_i\in\R^n  \mbox{ and } \sum_{i=0}^m\mu_i=1 \right\},
	\]
which lies in a $d$-dimensional affine subspace.	
If each $v_i\in\Z^n$ (resp. $v_i\in\Q^n$), we say that $\Pc$ is a \emph{lattice} polytope (resp. \emph{rational} polytope).	
Alternatively, $\Pc$ can be expressed	as the intersection of finitely many halfspaces in $\R^n$.

The \emph{lattice point enumerator} of $\Pc$, where $d=\dim(P)\leq n$, is the function 
	\[
	i(\Pc,t)=\# (t\cdot\Pc\cap \Z^n)
	\]
where $t\cdot \Pc=\{t\cdot \alpha \ : \ \alpha\in\Pc\}$ is the $t$th dilate of $\Pc$ with $t\in\Znn$.  By theorems of Ehrhart \cite{Ehrhart}, if $\Pc$ is lattice, $i(\Pc,t)$ is a polynomial in the variable $t$ of degree $d$ and  if $\Pc$ is rational, $i(\Pc,t)$ is a quasipolynomial in the variable $t$ of degree $d$. Subsequently, we will call $i(\Pc,t)$ the \emph{Ehrhart polynomial} of $\Pc$ or the \emph{Ehrhart quasipolynomial} of $\Pc$ in each respective case.
We say that $\Pc$ is \emph{Ehrhart positive} if each coefficient of $i(\Pc,t)$ is a positive rational number.
For background on Ehrhart positivity, see the survey of Liu \cite{Liu-Survey}.

We should note that given a rational or lattice polytope $\Pc$, we can associate a rational cone,  and thus an affine, normal semigroup algebra, to $\Pc$.
The \emph{cone over $\Pc$} is 
	\[
	\operatorname{cone}(\Pc):=\operatorname{span}_{\R_{\geq 0}}\left\{(\p,1) \, : \, \p\in\Pc\cap\Z^n \right\}
	\] 
and the \emph{affine, normal semigroup algebra} of $\Pc$ is $\C[\Pc]:=\C[\operatorname{cone}(\Pc)\cap \Z^{n+1}]$.
We can alternatively view $i(\Pc,t)$ as the \emph{Hiblert function} or $\C[\Pc]$.

Suppose that $\Pc$ is lattice. It is often more desirable to study the Ehrhart polynomial under a suitable change of basis. In particular, let $h^\ast(\Pc;z)$ be the Ehrhart polynomial expressed in the polynomial basis $\left\lbrace z^i:={{t+n-i} \choose n} \right\rbrace_{i=0}^n$.
The polynomial $h^\ast(\Pc;z)$ is called the \emph{$h^\ast$-polynomial} of $\Pc$.
Equivalently, we can view the  $h^\ast$-polynomial as the numerator of the \emph{Ehrhart series} of $\Pc$ which is
	\[
	\operatorname{Ehr}_\Pc(z):=\sum_{t\geq 0}i(\Pc,t)z^t=\frac{h^\ast(P;z)}{(1-z)^{d+1}}.
	\]

In the special case of a $d$-simplex $\Delta\subset\R^n$, we have a geometric interpretation of $h^\ast(\Delta;z)$. 
The \emph{(half-open) fundamental parallelpiped} of $\Delta$ is 
	\[
	\Pi_\Delta:=\left\lbrace \sum_{i=0}^d \mu_i(v_i,1)\in\R^{n+1} \, : \, 0\leq \mu_i<1 \, , \, i\in\{0,\ldots, d\} \right\rbrace.
	\]
We then have 
	\[
	h^\ast(\Delta;z)=\sum_{(x_1,\ldots, x_n,x_{n+1})\in \Pi_\Delta\cap \Z^{n+1}} z^{x_{n+1}}.
	\]		

We say that a lattice polytope $\Pc$ has the \emph{integer decomposition property} (IDP) if for any $\x\in t \Pc\cap \Z^n$, there exists $t$ lattice points $\{\p_1,\p_2,\ldots,\p_t\}\in P\cap \Z^n$ such that $\p_1+\p_2+\cdots +\p_t=\x$. 
In this case, we say that $\Pc$ has the IDP. 
This is equivalent to saying the semigroup algebra $\C[\Pc]$ is generated entirely in degree 1.

A \emph{triangulation} $\mathcal{T}$ of a lattice polytope $\Pc$ of dimension $d$ is a subdivision of $\Pc$ into $d$-dimensional lattice simplices framing a simplicial complex. 
We say that a triangulation $\mathcal{T}$  is \emph{regular} if each simplex  in $\mathcal{T}$ is a domains of linearity of a convex piecewise linear function.
We say that a triangulation $\mathcal{T}$ is \emph{unimodular} if each simplex in $\mathcal{T}$ is a unimodular simplex.
If $\Pc$ has a  unimodular triangulation, then $\Pc$ has the IDP.

Let $\Pc$ be a lattice polytope of dimension $d$ such that $\Pc\setminus \partial\Pc \cap \Z^n=\{{\bf 0}\}$.
We say that $\Pc$ is \emph{reflexive} if the dual polytope 
	\[
	\Pc^\vee:= \left\{y\in\R^n \, : \, 	\langle x, y\rangle \leq 1 \, \mbox{ for all } x\in \Pc \right\}
	\]
is a lattice polytope. 
It is a well known result of Hibi \cite{HibiDualPolytopes} that the following are equivalent:
	\begin{itemize}
	\item $\Pc$ is reflexive (up to translation)
	\item $h^\ast(\Pc;z)$ is a palindromic polynomial of degree $d$.
	\end{itemize}
	
Given $c\in\Zpos$, $\Pc$ is \emph{Gorenstein of index $c$} if the the polytope $c\Pc$ is reflexive. Equivalently, $\Pc$ is Gorenstein of index $c$ (up to unimodular equivalence) if 	$h^\ast (\Pc; z)$ is a palindromic polynomial of degree $d-c+1$ \cite{DeNegriHibi}. 
One should additionally, note that this precisely means that the  $\operatorname{cone}(\Pc)$ is a Gorenstein cone. 
That is, that there exists some $\c\in\Z^{n+1}$ such that $\c+(\cone(\Pc)\cap\Z^{n+1})=\cone(\Pc)^\circ \cap \Z^{n+1}$. Hence, $\Pc$ is Gorenstein if and only if the algebra $\C[\Pc]$ is Gorenstein.

A generalization of the Gorenstein property is the level property.
We say that $\Pc$ is \emph{level} if there exists some $\c_1,\ldots,\c_k\in \Z^{n+1}$ such that 
	\[
	\cone(\Pc)^\circ\cap \Z^{n+1}=\sum_{i=1}^{k}(\c_i+(\cone(\Pc)\cap \Z^{n+1})
	\] 
and $\c_{1_{n+1}}=\cdots=\c_{k_{n+1}}$.
Note that in the case of $k=1$, this is precisely the definition of Gorenstein. 
This implies that the algebra $\C[\Pc]$ is a \emph{level semigroup algebra}, which generalizes the notion of a Gorenstein algebra. 
While palidromicity of $h^\ast(\Pc;z)$ is lost with this generalization, additional inequalities on the coefficents remain. 
We recommend \cite{StanleyGreenBook} for details on level algebras.

A commonly studied property in Ehrhart theory is $h^\ast$-unimodality. That is, determining under what conditions the polynomial $h^\ast(\Pc;z)$ has a unimodal sequence of coefficients.
In general, one should not expect $h^\ast(\Pc;z)$ to be unimodal, even under the assumptions of Gorenstein (see, e.g., \cite{PayneNonUnimodal}). 
However, the following theorem on unimodality, due to Bruns and R\"omer, is known.

\begin{theorem}[\cite{BrunsRomer}]
If $\Pc$ is Gorenstein and has a regular unimodular triangulation, then $h^\ast(\Pc;z)$ is unimodal. 
\end{theorem} 

The following conjecture of Hibi and Ohsugi, which is a weakening of the Bruns and R\"omer result, is the subject of much research.

\begin{conjecture}[\cite{HibiOhsugi-Conjecture}]\label{conj:HibiOhsugi}
Let $\Pc$ be Gorenstein and have the IDP, then $h^\ast (\Pc; z)$ is unimodal.
\end{conjecture}

There is even a weaker question, which does not have a negative answer as of yet.

\begin{question}[\cite{SchepersVanLangenoven}]\label{ques:IDPUnimodal}
Does every IDP polytope have a unimodal $h^\ast$-polynomial?
\end{question}

\subsection{Local Ehrhart Theory}

Building on the classical $h^\ast$-polynomials in Ehrhart theory is the theory of \emph{local $h^\ast$-polynomials}. 
This theory has developed in large part due to Betke and McMullen \cite{BetkeMcMullen}, Katz and Stapledon \cite{KatzStapledon-Ehrhart}, Schepers and Van Langenhoven \cite{SchepersVanLangenoven}, and Stanley \cite{Stanley-localh}.

Suppose that $\Delta$ is a lattice $d$-simplex in $\R^n$ with vertex set $\{v_0,\ldots,v_d\}$. Then the \emph{open-parallelpiped} of $\Delta$ is
	\[
	\Pi^\circ_\Delta:=\left\lbrace \sum_{i=0}^d \mu_i(v_i,1)\in\R^{n+1} \, : \, 0<\mu_i<1 \, , \, i\in\{0,\ldots, d\} \right\rbrace.
	\]
Then the \emph{local $h^\ast$-polynomial} for the simplex $\Delta$ is
	\[
	\ell^\ast(\Delta;z):= \sum_{(x_1,\ldots,x_n,x_{n+1})\in {\Pi^\circ_\Delta}\cap \Z^{n+1}}z^{x_{n+1}}.
	\]	
We note1 that this is similar in nature to the geometric interpretation of $h^\ast(\Delta;z)$.
The polynomial $\ell^\ast(\Delta;z)$ is sometimes referred to as the  \emph{box polynomial} of $\Delta$. 
We should note that there is an involution on $\Pi^\circ_\Delta$ which sends elements at height $i$ to height $d-i$, which implies that $\ell^\ast(\Delta;z)$ is a palindromic polynomial.
The following theorem relates the local $h^\ast$-polynomial and the classic $h^\ast$-polynomial. This theorem uses the notion of an \emph{$h$-polynomial} of a simplicial complex. For background and details, please consult \cite{Stanley-EC1}.

\begin{theorem}[{\cite{BetkeMcMullen}}]
Let $\Pc$ be a lattice polytope with lattice triangulation $T$. Then
	\[
	h^\ast(\Pc;z)=\sum_{\Delta\in T}h(\operatorname{link}_T(\Delta);z)\ell^\ast(\Delta;z).
	\]
\end{theorem}
In the case of $\Pc$ a reflexive polytope, there is the following theorem. 
\begin{theorem}[\cite{BetkeMcMullen}]
Let $\Pc$ be a reflexive lattice polytope with lattice triangulation $T$ of the boundary $\partial \Pc$. Then
	\[
	h^\ast(\Pc;z)=\sum_{\Delta\in T}h(\operatorname{link}_T(\Delta);z)\ell^\ast(\Delta;z).
	\]
\end{theorem}

We say that a lattice triangulation $T$ of $\Pc$ is \emph{box unimodal} if it is regular and  $\ell^\ast(\Delta;z)$ is unimodal for $\Delta\in T$.
Therefore, if a reflexive polytope $\Pc$ admits a boundary triangulation whose simplices all have unimodal local $h^\ast$-polynomials, then $h^\ast(\Pc; z)$ is also unimodal.

In the case of $\Pc$ not a simplex, one can still define the local $h^\ast$-polynomial $\ell^\ast(\Pc;z)$ with some additional work. 
Let $F(\Pc)$ denote the set of faces of $\Pc$. Note that it is clear that $F(P)$ is a ranked poset under inclusion of faces where the rank is $\rho(F)=\dim(F)+1$.
The \emph{$g$-polynomial} of $\Pc$ is defined recursively as follows. If $\Pc=\emptyset$, then $g(\Pc;z)=1$. If $\Pc$ is dimension $d,$ then $g(\Pc; z)$ is the unique polynomial of degree strictly less than $\frac{d}{2}$ satisfying
	\[
	z^d g\left(\Pc;\frac{1}{z} \right)= \sum_{F\in F(\Pc)} g(F;z) (t-1)^{d-\dim(F)-1}.
	\]
The reader should consult \cite{Stanley-EC1} for details on the $g$-polynomial.
Then the \emph{local $h^\ast$-polynomial} of $\Pc$ is
	\[
	\ell^\ast(\Pc;z)=\sum_{F\in F(\Pc)}(-1)^{\dim(\Pc)-\dim(F)}h^\ast(F;z)g([F,\Pc]^\vee;z),
	\]
where $[F,\Pc]^\vee$ is polytope obtained from the dual of the interval $[F,\Pc]$ in $F(\Pc)$.	
\subsection{Polyhedral objects arising from lecture hall partitions}
We now briefly define the common polyhedral objects arising from lecture hall partitions.
We begin with the \emph{$\s$-lecture hall cone}.

\begin{definition}
Let $\s\in\Zpos^n$. The \emph{$\s$-lecture hall cone} is
	\[
	\Cc_n^{(\s)}:=\left\{ \lambda\in\R^n \, : \, 0\leq \frac{\lambda_1}{s_1}\leq \cdots\leq \frac{\lambda_n}{s_n} \right\},
	\]
which may equivalently be expressed as 
	\[
	\Cc_n^{(\s)}:=\operatorname{span}_{\R_{\geq 0}}\left\{ (0,\ldots,0, s_i, \ldots, s_n) \, : \, 1\leq i\leq n \right\}.
	\]	
\end{definition}
It is straightforward to see that $\Cc_n^{(\s)}\cap\Z^n=L_n^{(\s)}$. 
Subsequently, studying this polyhedral object yields the same data as studying the lecture hall partitions directly.
Results on lecture hall cones include a classification of the Gorenstein property, characterization of Hilbert bases in special case, and triangulation results for the cone $\Cc_n^{(1,2,\ldots,n)}$. All of these results will discussed in detail in Section \ref{sec:Cones}.

A related object are the \emph{$\s$-lecture hall simplices}\footnote{Much of the literature refers to these as \emph{lecture hall polytopes.} However, since these polytopes are simplices and since the development of lecture hall order polytopes, this is a better name in terms of specificity.}, which arise as a truncation of the lecture hall cone. Particularly, they can be defined as follows.

\begin{definition}
Let $\s\in\Zpos^n$. The \emph{$\s$-lecture hall simplex} is
	\[
	\Pbn:=\left\{ \lambda\in\R^n \, : \, 0\leq \frac{\lambda_1}{s_1}\leq \cdots\leq \frac{\lambda_n}{s_n}\leq 1 \right\},
	\]
which may equivalently expressed as 
	\[
	\Pbn:=\conv\left\{(0,\ldots,0), (0,\ldots,0,s_i,\ldots,s_n) \, : \, 1\leq i \leq n \right\}.
	\]	
\end{definition}
This family of simplices has been very well-studied. 
There are results regarding the $h^\ast$- polynomials, existence of the IDP, existence of triangulations, classification of the Gorenstein property and the level property, and results on the local $h^\ast$-polynomial. 
These results will be outlined in detail in Section \ref{sec:Simplices}. 

A less well-studied object arising from lecture hall partitions are the \emph{rational $\s$-lecture hall simplices}.
\begin{definition}
Let $\s\in\Zpos^n$. The \emph{rational $\s$-lecture hall simplex} is
	\[
	R_n^{(\s)}:=\left\{ \lambda\in\R^n \, : \, 0\leq \frac{\lambda_1}{s_1}\leq \cdots\leq \frac{\lambda_n}{s_n}\leq \frac{1}{s_n} \right\}.
	\]
\end{definition}
The study of these objects was initiated by Pensyl and Savage \cite{PensylSavage-Rational}.
This theory has been used while studying both the $\s$-lecture hall cones and $\s$-lecture hall simplices.
Given the relative lack of known results regarding these rational polytopes, there is no section specifically dedicated to them in this survey, but they appear in  Section \ref{sec:Cones} as well as Section \ref{sec:Open}.

Another important polyhedral object arrising from lecture hall partitions are \emph{$\s$-lecture hall order polytopes}. These objects were introduced by Br\"{a}nd\'{e}n and Leander \cite{BrandenLeander-PPartition} as a generalization of the classical theory of \emph{order polytopes} defined by Stanley \cite{Stanley-PosetPolytopes} which have been well studied.

In order to define $\s$-lecture hall order polytopes, we first need the notion of an \emph{$\s$-lecture hall $P$-parition.}
Let $P=([n],\preceq)$ be a labeled poset on $[n]$.
We say that $P$ is \emph{naturally-labeled} if $i\preceq j$ in $P$ implies that $i\leq j$ in $\Z$. 

\begin{definition}
Let $P=([n],\preceq)$ be a labeled poset  and let $\s:[n]\to \Zpos$ be an arbitrary map. A \emph{lecture hall $P$-partition} is a map $f:[n]\to \R$ such that 
	\begin{enumerate}
	\item if $x \prec y$, then $f(x)/s(x) \leq f(y)/s(y)$, and 
	\item if $x\prec y$ and $x>y$, the $f(x)/s(x) < f(y)/s(y)$.	
	\end{enumerate}			  
\end{definition}
If we consider naturally-labeled posets, we need only consider the first condition of the definition. 
We now define the lecture hall order polytopes.

\begin{definition}\label{def:sorderpoly}
Let $P$ be a naturally-labeled poset. The \emph{$\s$-lecture hall order polytope} associated to $(P,\s)$ is 
	\[
	O(P,\s):=\{f\in\R^n \, : \, f \mbox{ is a } (P,\s)\mbox{-partition and } 0\leq f(x)/s(x)\leq 1 \mbox{ for all } x\in [n]\} 
	\]
\end{definition}
\begin{remark}
Note that the naturally-labeled condition in Definition \ref{def:sorderpoly} ensures that the resulting polytopes are closed.
One could define a similar notion for any lecture hall $P$-partition, however, the resulting polytope will be partially open without the naturally-labeled condition (i.e. missing certain facets).
\end{remark}
We should note that if we let $P$ be a totally ordered chain, we recover the $\s$-lecture hall simplex as $O(P,\s)=\Pbn$.
While these objects are relatively recent in development, significant results have been obtained. 
Specifically, there are known results about the $h^\ast$ polynomial including some results real-rootedness. 
There are also known results regarding box unimodal traingulations and local $h^\ast$-polynomials.
These results will be outlined in Section \ref{sec:order}.

\section{Results:  Lecture hall cones}\label{sec:Cones}

In this section, we discuss the results on $\s$-lecture hall cones. 
The known results fall into three categories: the classification of the Gorenstein property, the characterization of Hilbert bases for certain Gorenstein lecture hall cones, and the special case of the lecture hall cone $\Cc_n^{(1,2,\ldots,n)}$.

\subsection{Gorenstein characterization}

In combinatorial commutative algebra and polyhedral geometry, Gorenstein objects are well-behaved and it is desirable to classify which objects in the family have this property. 
In the case of $\s$-lecture hall cones, this property is entirely controlled by number theoretic properties of the $\s$-sequence.  

\begin{theorem}[{\cite{BeckEtAl-GorensteinLHC,BME-LHP2}}] \label{GorensteinCriteria}
For a positive integer sequence $\s$, the $\s$-lecture hall  cone $\Cc_n^{(\s)}$ is Gorenstein if and only if there exists some $\c \in \Z^n$ satisfying
	\[
	c_j s_{j-1}=c_{j-1}s_j+\gcd(s_j,s_{j+1})
	\]
for $j>1$, with $c_1=1$. 	
\end{theorem}
This  result is implicit from the seminal work of Bousquet-M\'{e}lou and Eriksson \cite{BME-LHP2}, but was formulated in this context by Beck et al. \cite{BeckEtAl-GorensteinLHC}.
While this completely classifies when lecture hall cones are Gorenstein, the implicit nature of the definition is somewhat mysterious. 
Subsequently, it is of interest to study particular types of sequences, as the specificity may provide a more enlightening answer.
Given a sequence $\s\in\Zpos^n$, we say that $\s$ is \emph{$\u$-generated} if there is a sequence $\u=(u_1,u_2,\cdots,u_{n-1})\in\Zpos^{n-1}$ such that $s_2=u_1s_1-1$ and $s_{i+1}=u_is_i-s_{i-1}$ for $i>1$.
Note that by construction a $\u$-generated sequence has the property $\gcd(s_i,s_{i+1})=1$ for $1\leq i< n$.
We can now state a characterization of the Gorenstein property for the special case where $\gcd(s_i,s_{i+1})=1$ for all $i$.

\begin{theorem}[{\cite{BeckEtAl-GorensteinLHC,BME-LHP2}}]\label{PairwiseGorenstein}
Let $\s=(s_1,\cdots,s_n)$ be a sequence of positive integers  such that $\gcd(s_i,s_{i+1})=1$ for $1\leq i<n$. Then $\Cc_n^{(\s)}$ is Gorenstein if and only if $\s$ is $\u$-generated by some sequence $\u=(u_1,u_2,\cdots,u_{n-1})$ of positive integers. When such a sequence exists, the Gorenstein point $\c$ for $\Cc_n^{(\s)}$ is defined by $c_1=1$, $c_2=u_1$, and for $2\leq i<n$, $c_{i+1}=u_ic_i-c_{i-1}$.
\end{theorem}

Additionally, the two-term recurrence case has been studied. 
The following theorem specifies precisely when such a sequence will produce a Gorenstein lecture hall cone.

\begin{theorem}[{\cite{BeckEtAl-GorensteinLHC}}]
Let $\ell>0$ and $b\neq 0$ be integers such that $\ell^2+4b\geq 0$. Let $\s=(s_1,s_2,\ldots)$ defined by 
	\[
	s_n=\ell s_{n-1}+bs_{n-2}
	\]
with $s_1=1$ and $s_0=0$. The $\Cc_n^{(\s)}$ is Gorenstein for all $n\geq 0$ if and only if $b=-1$
If $b\neq -1$, there exists $n_0$ such that $\Cc_n^{(\s)}$ fails to be Gorenstein for all $n\geq n_0$. 	
\end{theorem}

\subsection{Hilbert basis results}
Given a rational polyhedral cone, determining the Hilbert basis is a natural line of study. 
Understanding the Hilbert bases for $\s$-lecture hall cones is particularly of interest, because this is equivalent to determine a minimum additive generating set for the set of all $\s$-lecture hall partitions.
In the special case of $\s=(1,2,\ldots,n)$, the following elegant description of the Hilbert basis is known.

\begin{theorem}[{\cite{BeckEtAl-TriangulationsLHC}}]\label{BBKSZHilbert}
For each $A=\{a_1<a_2<\cdots<a_k\}\subseteq[n-1]$, define the element $v_A$ to be
	\[
	v_A=(0,\ldots,0, a_1,a_2,\ldots, a_k, a_k+1).
	\]	
The Hilbert basis for $L_n^{(1,2,\ldots,n)}$ is 
	\[
	\Hc_n^{(1,2,\ldots,n)}:=\{v_A \ : \ A \subseteq [n-1]\}.
	\]
As a corollary, the semigroup algebra $\C[\Cc_n^{(1,2,\ldots,n)}]$ is generated entirely by elements in degree $1$ with respect to the grading given by $\lambda\mapsto(\lambda_n-\lambda_{n-1})$. 
\end{theorem}

Given this surprisingly nice, combinatorial description for the Hilbert basis in this case, a natural follow-up is to determine which other $\s$-sequences have combinatorially nice Hilbert bases.
It is unlikely that one could give a  universal description for the Hilbert basis of any $\s$-lecture hall cone. 
Even in the case of 2-dimensional cones, one cannot even bound the cardinality of the Hilbert basis beyond the trivial bounds for simplicial cones.
Subsequently, to hope for meaningful results, it is necessary to restrict our classes of $\s$-sequences.

The first more general family of $\s$-sequences considered in the literature on Hilbert bases are the \emph{$1\bmod k$-sequences}.
Given a $k\in\Zpos$, the $1\bmod k$-sequence of length $n$ is
	\[
	\s=(1,k+1,2k+1,\ldots, (n-1)k+1).
	\]
For simplicity of notation, when $\s$ is as above, we will use $L_{k,n}$, $\Cc_{k,n}$, and $\Pb_{k,n}$ to denote the lecture hall partitions, lecture hall cone, and lecture hall simplex respectively.
We have the following result for the Hilbert basis of $\Cc_{k,n}$.	

\begin{theorem}[{\cite{Olsen-HilbertBases}}]
For all $k\geq 1$, the Hilbert basis for the $1\,\operatorname{mod}\, k$ cones in $\R^n$, denoted $\Cc_{k,n}$, consist of the following elements:
	\begin{itemize}
	\item The element $v_A:=(0,0,\ldots,0, a_1,a_2,\ldots,a_k,a_k+1)$ for each $A\subseteq[n-2]$ where\\  $A=\{a_1<a_2<\cdots<a_k\}$;
	\item Element $w\in L_{k,n}$, where $w_{n-1}=(n-2)k+1$ and $w_n=(n-1)k+1$;  
	\end{itemize}
where $L_{n,k}$ denotes the set of $1\,\operatorname{mod}\, k$ lecture hall partitions. 
\end{theorem}
When $k=1$, we recover the sequence $\s=(1,2,\ldots,n)$ and moreover, this result agrees with previous characterization of Beck et al. \cite{BeckEtAl-TriangulationsLHC}.
We can also enumerate the cardinality of the Hilbert basis of this cone with the following algebraic expression.

\begin{corollary}[{\cite{Olsen-HilbertBases}}]
	\[
	|\Hc_{k,n}|=\frac{(k+1)^{n-2}+(k-1)}{k}+2^{n-2}.
	\]
\end{corollary}
To see this result, it is clear that the $2^{n-2}$ appears from enumerating the elements indexed by subsets.
The remaining elements can be enumerated by observing that each Hilbert element can be realized as a lattice point in $\Pb_{k,n-2}$ and applying Ehrhart theory results of Savage and Visawanathan \cite{SavageViswanathan-1/kEulerian}.

The second more general family of $\s$-sequences studied are the \emph{$\ell$-sequences}. 
These sequences are recursively defined by $s_1=1$, $s_2=\ell$, and $s_i=\ell s_{i-1}-s_{i-2}$ for $i\geq 2$.
For ease of notation, given the above sequence let $L_n^\ell$, $\Cc_n^\ell$, $\Pb_n^\ell$, and $R_n^\ell$ denote the lecture hall partitions, lecture hall cone, lecture hall simplex, and rational lecture hall simplex respectively.
The Hilbert basis results for these lecture hall cones are given as follows. 

\begin{theorem}[{\cite{Olsen-HilbertBases}}]
Let $\s=(s_1,s_2,\ldots, s_n)$ be an $\ell$-sequence for some $\ell\geq 2$. The Hilbert basis $\Hc_n^\ell$ for the $\ell$-sequence cone  $\Cc_n^\ell$ is 

	 \[
	 \Hc_n^\ell=\bigcup_{i=0}^{n}\left\{\lambda\in L_n^\ell \ : \ \lambda_{n-1}=s_i \, , \, \lambda_n=s_{i+1}\right\}
	 \]
where $L_n^\ell$ denotes the set of $\ell$-sequence lecture hall partitions.  	 
	\end{theorem}
We should note that when $\ell=2$, we recover the sequence $\s=(1,2,\ldots,n)$ and this characterization agrees with the Beck et al. result \cite{BeckEtAl-TriangulationsLHC}.
By identifying these Hilbert basis elements with lattice points in dilates of rational lecture hall simplices, one can derive the following expression for the cardinality of the Hilbert basis.	

\begin{corollary}[{\cite[Corollary 4.2]{Olsen-HilbertBases}}]
	\[
	|\Hc_n^\ell|=2+\sum_{j=1}^{n-2}i(R_{n-2}^\ell,s_j)
	\] 
where $i(R_{n-2}^\ell,t)$ denotes the Ehrhart quasipolynomial of the rational lecture hall simplices $R_{n-2}^\ell$.	 
\end{corollary}	

In addition to these Hilbert basis results, the author has also classified the Hilbert bases for $\u$-generated Gorenstein lecture hall cones when $n=2,3,4$ \cite{Olsen-HilbertBases}.

\subsection{Triangulation of $L_n^{(1,2,\ldots,n)}$}
In the specific scenario of $\s=(1,2,\ldots,n)$, one can make more specific observation of the geometric structure of the lecture hall cone. 
The result of Theorem \ref{BBKSZHilbert} verifies that the Hilbert basis of $\Cc_n^{(1,2,\ldots,n)}$ is entirely at height 1 with respect to a chosen grading. 
This motivates asking whether or not that cone has a regular, unimodular triangulation.
Beck et al. verify that  this cone can be represented as the cone over a polytope. 
Specifically, $\Cc_n^{(1,2,\ldots,n)}\cong\cone(\mathscr{R}_n)$, where 
	\[
	\mathscr{R}_n=\conv
	\begin{bmatrix}
	1 & 1 & 1 & \cdots & 1\\
	0 & n-1 & 1 & \cdots & 1\\
	0 & 0 & n-2 & \cdots & 1\\
	\vdots & \vdots & & \ddots & \vdots\\
	0 & 0 & 0 & \cdots & 1
	\end{bmatrix}.
	\]
They then state the following theorem.
\begin{theorem}[\cite{BeckEtAl-TriangulationsLHC}]
For $n\geq 1$, $\mathscr{R}_n$ has a regular, flag, unimodular triangulation. 
Thus $\Cc_n^{(1,2,\ldots,n)}$ has a regular, unimodular triangulation.
\end{theorem}

\section{Results: Lecture hall simplices}\label{sec:Simplices}
We now discuss many of the known results for a lecture hall simplex $\Pbn$. 
Perhaps the most pervasive question in Ehrhart theory is: given in a lattice polytope $\Pc$, is there a combinatorial interpretation for $h^\ast(\Pc;z)$?
To determine such an interpretation for $\Pbn$, we must define some statistics related to lecture hall partitions.
The set of \emph{$\s$-inversion sequences} is
	\[
	\Ibn:=\left\lbrace\e\in\Znn^n \, : \, 0\leq e_i<s_i \mbox{ for } 1\leq i\leq n \right\rbrace.
	\]
Given $\e\in\Ibn$, the \emph{ascent set} of $\e$ is 
	\[
	\Asc(\e):=\left\lbrace i\in\{0,1,\ldots,n-1\} \, : \frac{e_i}{s_i}< \frac{e_{i+1}}{s_{i+1}} \right\rbrace
	\]
where by convention $s_0=1$ and $e_0=0$.
The \emph{ascent number} of $\e\in\Ibn$ is $\asc(\e):=\#\Asc(\e)$.		
Savage and Schuster show the following result.

\begin{theorem} [{\cite{SavageSchuster}}]
For any $\s\in\Zpos^n$, 
	\[
	h^\ast\left(\Pbn;z \right)=\sum_{\e\in\Ibn}z^{\asc(\e)}.
	\]
\end{theorem}
The proof of this result is essentially a bijection between lattice points in $\Pi_{\Pbn}$ and $\s$-inversion sequences with respect to four statistics.
When restricting to the ascent statistic alone, this corresponds to height of lattice points. 
In the special case of $\s=(1,2,\ldots,n)$, one can show that 
	\[
	h^\ast(\Pb_n^{(1,\ldots,n)};z)=\sum_{\e\in\Ibn}z^{\asc(\e)}=\sum_{\pi\in \mathfrak{S}_n}z^{\des(\pi)}=A_n(z)
	\]
where $\mathfrak{S}_n$ is the symmetric group and thus $A_n(z)$ is the \emph{Eulerian Polynomial} (see, e.g., \cite{Stanley-EC1} for more details on Eulerian polynomials).
Subsequently, $h^\ast(\Pbn;z)$ generalize this family and are known as the \emph{$\s$-Eulerian polynomials}.

The $\s$-Eulerian polynomials also possess many of the same nice properties of Eulerian polynomials, namely unimodality and real-rootedness due to the following theorem of Savage and Visontai \cite{SavageVisontai}. 

\begin{theorem}[{\cite{SavageVisontai}}]
\label{thm:sEulerRealRoot}
For any $\s\in\Zpos$, the $\s$-Eulerian polynomial $h^\ast(\Pbn;z)$ has only real roots. Hence, the coefficients of $\s$-Eulerian polynomials are log-concave and unimodal. 
\end{theorem}
The proof of this theorem uses the idea of compatible polynomials and interlacing to prove real-rootedness. 
For background on real-rootedness in a combinatorial setting, the reader should consult \cite{Branden-Survey}. 
Unfortunately, $\s$-Eulerian polynomials do not apriori possess the same palindromic properties as the Eulerian polynomial $A_n(z)$.
Classification of this property is equivalent to the classification of which lecture hall simplices $\Pbn$ are Gorenstein, which will be discussed in a later subsection.

While $h^\ast(\Pbn;z)$ is unimodal for any choice of $\s$, it is not true in general that $\Pbn$ is Ehrhart positive.
For $a,b,k_2\in\Zpos$ and $k_1,k_3\in\Znn$, let the vector $(1^{k_1},a,1^{k_2},b,1^{k_3})=(\underbrace{1,\ldots,1}_{k_1},a,\underbrace{1,\ldots,1}_{k_2},b,\underbrace{1,\ldots,1}_{k_3})$. 
The following theorem of Liu and Solus \cite{LiuSolus} provide examples where the Ehrhart polynomial of $\Pbn$ has negative coefficients.

\begin{theorem}[\cite{LiuSolus}]\label{thm:EhrhartPos}
For every $n\ge 3$, there exist $\s$-lecture hall simplices of the form $\Pb_n^{ (1^{k_1},a,1^{k_2},b,1^{k_3})}$ which are not Ehrhart positive.
\end{theorem}

In addition to understanding this classical Ehrhart data of $\Pbn$, it is also of interest to investigate geometric properties such as the existence of the integer decomposition property (IDP) and the existence of particular triangulations.
Regarding the integer decomposition property, Hibi Tsuchiya, and the author showed some partial results \cite{HibiOlsenTsuchiya-LHPGorensteinIDP}, specifically that $\Pbn$ has the IDP when $\s$ is a weakly monotonic sequence and that $\Pb_n^{(\s,1,\t)}$ has the IDP if and only if $\Pbn$ and $\Pb_n^{(\t)}$ are both have the IDP.
This results has since been generalized to all $\Pbn$ in forthcoming work of Br\"and\'{e}n and Solus \cite{BrandenSolus}, where they prove stronger results for the $\s$-lecture hall order polytopes.

\begin{theorem}[\cite{BrandenSolus}]
Let $\s\in\Zpos^n$ be a sequence. The $\s$-lecture hall simplex $\Pbn$ has the IDP.
\end{theorem}

Given that all lecture hall simplices possess the IDP, a natural extension would be the classification of which lecture hall simplices possess regular, unimodular triangulations.
There are two notable results in this direction which follow.

\begin{theorem}[{\cite{BrandenSolus}}]
Suppose that $\s$ is weakly monotone sequence such that $|s_{i+1}-s_{i}|\leq 1$. 
Then $\Pbn$ has a regular, unimodular triangulation.
\end{theorem}

\begin{theorem}[{\cite{HibiOlsenTsuchiya-LHPGorensteinIDP}}]
Suppose that $\Pbn$ has a regular, unimodular triangulation. Then $\Pb_{n+1}^{(\s,ks_n)}$ has a regular, unimodular triangulation for any $k\in\Zpos$.
\end{theorem}
The first theorem above, due to Br\"and\'{e}n and Solus \cite{BrandenSolus}, uses the theory of Gr\"obner bases. 
The second theorem, due to Hibi, Tsuchiya, and the author \cite{HibiOlsenTsuchiya-LHPGorensteinIDP}, uses the notion of a \emph{chimney polytope}.
The reader should consult \cite{HaaseEtAl-Triangulations} for details on triangulations, including these two methods.
There are no known examples of lecture hall simplices which do not admit a regular, unimodular triangulation.

\subsection{Gorenstein and level classifications}

Both the Gorenstein property and the level property for lecture hall simplices have been studied and a classification for both exists. 
The Gorenstein property was studied in limited cases by Hibi Tuschiya, and the author in \cite{HibiOlsenTsuchiya-LHPGorensteinIDP}.
These results have been fully generalized by Kohl and the author in \cite{KohlOlsen-Level}.

\begin{theorem}[{\cite{KohlOlsen-Level}}]
\label{thm:NewGorensteinChar}
Let $\s=(s_1,s_2,\ldots,s_n)\in \Z^{n}_{\geq 1}$. Then $\Pbn$ is Gorenstein if and only if there exists a $\c\in\Z^{n+1}$ satisfying
	\[
	c_js_{j-1}=c_{j-1}s_j+\gcd(s_{j-1},s_{j})
	\]
for $j>1$  and 
\[
c_{n+1}s_n = 1 + c_n
\]
with $c_1=1$.
\end{theorem}
While the next result is not as general, it guarantees that, under the condition that $\gcd(s_{i-1},s_i)=1$ for some $1<i \leq n$, the vertex cones of $\Pbn$ at $(0,0,\dots,0)$ and at $(s_1,s_2,\dots, s_n)$ being Gorenstein already implies that $\Pbn$ is Gorenstein.	 

\begin{theorem}[{\cite{KohlOlsen-Level}}]
\label{GorensteinChar}
Let $\s=(s_1,s_2,\ldots,s_n)$ be a sequence such that there exists some $1<i\leq n$ such that $\gcd(s_{i-1},s_i)=1$ and define $\overleftarrow{\s}=(\overleftarrow{s_1},\ldots,\overleftarrow{s_n}):=(s_n,s_{n-1},\ldots,s_1)$. Then $\Pbn$ is Gorenstein if and only if there exists $\c,\d\in\Z^n$ satisfying
	\[
	c_js_{j-1}=c_{j-1}s_j+\gcd(s_{j-1},s_{j})
	\]
and 
	\[
	d_j\overleftarrow{s_{j-1}}=d_{j-1}\overleftarrow{s_j}+\gcd(\overleftarrow{s_{j-1}},\overleftarrow{s_{j}})
	\]	
for $j>1$ with $c_1=d_1=1$.	
\end{theorem}
Note that Theorem \ref{thm:NewGorensteinChar} and Theorem \ref{GorensteinChar} provide a complete description for when $\s$-Eulerian polynomials are palindromic. See Table \ref{table:s-Eulerian} for some examples of palindromic $\s$-Eulerian polynomials which fall under the purview specifically of Theorem \ref{GorensteinChar}. 
Additionally, it is worth noting that in the special case of $\s=(2,3,\ldots,n+1)$, Hibi, Tsuchiya, and the author \cite{HibiOlsenTsuchiya-SelfDual} show that $\Pb_n^{(2,3,\ldots,n+1)}$ is a \emph{self-dual} reflexive simplex. That is, $\Pb_n^{(2,3,\ldots,n+1)}\cong{\Pb_n^{(2,3,\ldots,n+1)}}^\vee$.  
Moreover, $h^\ast(\Pb_n^{(2,3,\ldots,n+1)};z)$ is the usual Eulerian polynomial $A_{n+1}(z)$.

\begin{table}[h!]\label{table:s-Eulerian}
\centering
\footnotesize
	\begin{tabular}{|c|c|c|c|c|}
	\hline
	 & sequence $\s$ & corresponding $\c$ & corresponding $\d$ & $\s$-Eulerian polynomial\\
	 \hline\hline
	  & & & &  \\
	 (i) & $(2,1,3,2,1)$ & $(1,1,4,3,2)$ & $(1,3,5,2,5)$ & $1+5z+5z^2+z^3$\\
	 & & &  & \\
	 \hline
	  & & & & \\
	 (ii) & $(3,2,3,1,2)$ & $(1,1,2,1,3)$ & $(1,1,4,3,5)$ & $1+9z+16z^2+9z^3+z^4$\\ 
	 & & & & \\
	 \hline
	 & & & & \\
	 (iii) & $(1,4,3,2,3)$ & $(1,5,4,3,5)$ & $(1,1,2,3,1)$ & $1+ 16z+38z^2+16z^3+z^4$\\
	 & & & & \\
	 \hline
	 & & & & \\
	 (iv) & $(3,5,2,3,1)$ & $(1,2,1,2,1)$ & $(1,4,3,8,5)$ & $1+20z+48z^2+20z^3+z^4$\\
	 & & & & \\
	 \hline
	 & & & & \\
	 (v) & $(1,2,3,4,5)$ & $(1,3,5,7,9)$ & $(1,1,1,1,1)$ & $1+26z+66z^2+26z^3+z^4$\\
	 & & & &\\
	 \hline
	 & & & & \\
	 (vi) & $(1,2,5,8,3)$ & $(1,3,8,13,5)$ & $(1,3,2,1,1)$ & $1+50z+138z^2+50z^3+z^4$\\
	 & & & &  \\
	 \hline
	 & & & & \\
	 (vii) & $(4,3,2,5,3)$ & $(1,1,1,3,2)$ & $(1,2,1,2,3)$ &  $1+30z+149z^2+149z^3+30z^4+z^5$\\
	 & & & &\\
	 \hline
	 & & & & \\
	 (viii) & $(4,7,3,2,3)$ & $(1,2,1,1,2)$ & $(1,1,2,5,3)$ &  $1+43z+208z^2+208z^3+43z^4+z^5$\\
	 & & & & \\
	 \hline
	 & & & &  \\ 
	 (ix) & $(5,9,4,3,2)$ & $(1,2,1,1,1)$ & $(1,2,3,7,6)$ & $1+82z+457z^2+457z^3+82z^4+z^5$\\
	 & & & & \\
	 \hline
	 & & & & \\
	 (x) & $(3,5,12,7,2)$ & $(1,2,5,3,1)$ & $(1,4,7,3,2)$ & $1+175z+1084z^2+1084z^3+175z^4+z^5$\\
	 & & & & \\ 
	 \hline
	 & & & & \\
	 (xi) & $(3,11,8,5,2)$ & $(1,4,3,2,1)$ & $(1,3,5,7,2)$ & $1+180z+1139z^2+1139z^3+180z^4+z^5$\\
	 & & & & \\
	 \hline 
	 & & & &\\
	 (xii) & $(2,7,5,10,4)$ & $(1,4,3,7,3)$ & $(1,3,2,3,1)$ & $1+181z+1218z^2+1218z^3+181z^4+z^5$ \\
	 & & & &\\
	 \hline
	 & & & & \\
	 (xiii) & $(3,8,13,5,2)$ & $(1,3,5,2,1)$ & $(1,3,8,5,2)$ & $1+213z+1346z^2+1346z^3+213z^4+z^5$\\
	 & & & & \\
	 \hline
	\end{tabular}
	\medskip
\normalsize
\caption{Palindromic $\s$-Eulerian Polynomials  \cite{KohlOlsen-Level}}
\label{table:s-Eulerian}
\end{table}

In addition, Kohl and the author \cite{KohlOlsen-Level} are able to also provide a characterization of the level property in terms of $\s$-inversion sequences.

\begin{theorem}[{\cite{KohlOlsen-Level}}]
\label{LevelChar}
\commentout{[See Theorem \ref{LevelChar}]}
Let $\s=(s_1,s_2,\ldots, s_n)$ and let $r=\max\{\asc(\e) \, : \, \e\in \Ibn \}$. Then $\Pbn$ is level if and only if for any $\e\in\Ib_{n,k}^{(\s)}$ with $1\leq k <r$ there exists some $\e'\in \Ib_{n,1}^{(\s)}$ such that $(\e+\e')\in \Ib_{n,k+1}^{(\s)}$.
\end{theorem}
The proof of this theorem relies heavily on the fact that $\Pbn$ is a simplex and studying the additive properties of the lattice points in $\Pi_{\Pbn}$ using the bijection that the lattice points at height $i$ correspond to inversion sequences with $i$ inversions. 
The authors use this classification to prove that certain classes of lecture hall simplices are level and that we can construct level lecture hall simplicies of arbitrarily high dimension.
Additionally, levelness introduces new inequalities on the coefficients of the $\s$-Eulerian polynomials, as represented in following corollary.

\begin{corollary}[{\cite{KohlOlsen-Level}}]
Let $\s=(s_1,s_2,\ldots,s_n)$ be a sequence such that $\Pbn$ is level. Then the coefficients of the $\s$-Eulerian polynomial $h^\ast(\Pbn;z)=1+h_1^\ast z+\cdots+h_r^\ast z^r$ satisfy the the inequalities $h_i^\ast\leq h_j^\ast h_{i+j}^\ast$ for all pairs $i$ and $j$ such that $h_{i+j}^\ast>0$.
\end{corollary}

\subsection{Local $h^\ast$-polynomials for lecture hall simplices}

Given the breadth of knowledge on the traditional Ehrhart theory of $\s$-lecture hall simplices, a natural extension is to consider the local $h^\ast$-polynomials. Gustafsson and Solus in \cite{GustafssonSolus} do this in generality. 
To introduce these results, we need the following notation. Let 
	\[
	\tilde{I}_n^{(\s)}:=\left\lbrace (e_0,\ldots,e_{n+1})\in\Z^{n+2} \, : \, 0\leq e_i<s_i \mbox{ and } \frac{e_i}{s_i}\neq\frac{e_{i+1}}{s_{i+1}} \mbox{ for all } i\in\{0,\ldots,n\} \right\rbrace
	\]
where $s_0=s_{n+1}=1$ and $e_0=e_{n+1}=0$.	
	
\begin{theorem}[{\cite[Proposition 3.3]{GustafssonSolus}}]
Let $\s\in\Zpos^n$. The local $h^\ast$-polynomials for the $\s$-lecture hall simplex $\Pbn$ is 
	\[
	\ell^\ast(\Pbn;z)= d_n^\s(z):=\sum_{\e\in\tilde{I}_n^{(\s)}}z^{\asc(\e)}.
	\]		
\end{theorem}
The polynomials $d_n^\s(z)$ are called the \emph{$\s$-derangement polynomials}, because in the case of $\s=(2,3,\ldots,n+1)$, this polynomial agrees with the Eulerian polynomial computed only over $\pi\in\mathfrak{S}_n$ such that $\pi$ is a derangement.
Moreover, there is a parallel to the real-rootedness result of Savage and Visontai as follows.

\begin{theorem}[{\cite{GustafssonSolus}}]
Let $\s\in\Zpos$. Then the $\s$-derangement polynomial $d_n^\s(z)$ is real-rooted and hence log-concave and unimodal.
\end{theorem}

Suppose that $p(z)$ is a palindromic polynomial of degree at most $d$, then we can express $p$ in the basis $\{x^i(x+1)^{d-2i}\}_{i=0}^{\lfloor\frac{d}{2}\rfloor}$. That is, 
	\[
	p(z)=\sum_{i=0}^{\lfloor\frac{d}{2}\rfloor}\gamma_i x^i(x+1)^{d-2i}.
	\]
We say that $p(z)$ is \emph{$\gamma$-nonegative} if $\gamma_0,\ldots,\gamma_{\lfloor\frac{d}{2}\rfloor}$ are all nonnegative.
It is a straight forward observation to see that if $p(z)$ is $\gamma$-nonnegative, then $p(z)$ is unimodal.
The $\s$-derangement polynomials can be shown to have this property by the following theorem.

\begin{theorem}[{\cite{GustafssonSolus}}]
All $\s$-derangement polynomials are $\gamma$-nonnegatve.
\end{theorem}

\section{Results:  Lecture hall order polytopes}\label{sec:order}

The study of lecture hall order polytopes was initiated by Br\"and\'{e}n and Leander \cite{BrandenLeander-PPartition} by way of introducing \emph{lecture hall $P$-partitions}, which generalizes the $P$-partition theory of Stanley.
This subject has since been studied by Br\"and\'{e}n and Solus \cite{BrandenSolus} and Gustafsson and Solus \cite{GustafssonSolus} from a polyhedral geometry perspective, as well as by Corteel and Kim \cite{CorteelKim-Tableaux} from a partition theoretic framework. 
This section is organized as follows. First, we discuss general results regarding the integer decomposition property and the $h^\ast$-polynomials of $\s$-lecture hall order polytopes.
We then discuss results regarding unimodality and real-rootedness of $h^\ast$-polynomials. We conclude by discussing local Ehrhart theory results on these polytopes. 

We begin by discussing the results on the integer decomposition property.
The following theorem, due to Br\"and\'{e}n and Solus, further implies that both order polytopes and $\s$-lecture hall simplices have the IDP. 

\begin{theorem}[\cite{BrandenSolus}]\label{theorem:IDPOrder}
Let $P=([n],\preceq)$ be a labeled poset and $\s:[n]\to \Zpos$ be any function. Then $O(P,\s)$ has the IDP.
\end{theorem}

To introduce some of the initial Ehrhart results on $\s$-lecture hall order polytopes, we must review some theory of $P$-partitions and some generalizations. For additional background on $P$-partitions beyond that reviewed here, the reader should consult \cite{Stanley-EC1}.
Let $P=([n],\preceq)$ be a  labeled poset. Then the set 

	\[
	\Lc(P)=\left\lbrace \pi\in\Sn \, : \, \mbox{ if } \pi_i \preceq \pi_j, \mbox{ then }  i\leq j, \mbox{ for all } i,j\in[n] \right\rbrace
	\]
are the \emph{linear extensions} of $P$ or the \emph{Jordan-H\"older set} of $P$.
If $O(P)$ denotes the \emph{order polytope} of $P$, then the classical result of Stanley \cite{Stanley-PosetPolytopes} says that 
	\[
	h^\ast (O(P);z)=\sum_{\pi\in\Lc(P)}z^{\des(\pi)}.
	\] 
A similar result holds in the case of $\s$-lecture hall order polytopes. Given $\s:[n] \to \Zpos$, an \emph{$\s$-colored permutation} is a pair $\tau=(\pi,\r)$, where $\pi\in\Sn$ and $\r:[n]\to \Znn$ such that $\r(\pi_i)\in\{0,1,\ldots,\s(\pi_i)-1\}$ for all $1\leq i\leq n$.
Now, given a labeled poset $P$, 
	\[
	\Lc(P,\s):= \left\lbrace\tau \, : \, \tau=(\pi,\r) \mbox{ such that } \pi\in\Sn \mbox{ and } \tau \mbox{ is an } \s\mbox{-colored permutation} \right\rbrace.
	\]
Given $\tau=(\pi,\r)$, we say that $i\in[n-1]$ is a \emph{descent} of $\tau$ if 
	\[
	\begin{cases}
	\pi_i<\pi_{i+1} \mbox{ and } \r(\pi_i)/\s(\pi_i)>\r(\pi_{i+1})/\s(\pi_{i+1}), \mbox{ or },\\
	\pi_i>\pi_{i+1} \mbox{ and } \r(\pi_i)/\s(\pi_i)\geq \r(\pi_{i+1})/\s(\pi_{i+1})
	\end{cases}
	\] 
and we let 
	\[
	D_1(\tau):=\{i\in[n-1] \, : \, i \mbox{ is a descent}\}.
	\]
Moreover, let 
	\[
	D(\tau):=\begin{cases}
	D_1(\tau), & \mbox{ if } \r(\pi_n)= 0\\
	D_1(\tau)\cup\{n\}, & \mbox { if } \r (\pi_n)\neq 0 
	\end{cases}
	\]	
Given $\tau\in\Lc(P,\s)$, let $\des_\s(\tau)=|D(\tau)|$.
We can now state the following Ehrhart series result for $\s$-lecture hall order polytopes.  	

\begin{theorem}[{\cite[Corollary 3.7]{BrandenLeander-PPartition}}]
Let $P$ be a labeled poset and $\s:[n]\to \Zpos$. Then 
	\[
	h^\ast(O(P,\s);z)=\sum_{\tau\in \mathcal{L}(P,\s)} z^{\des_s(\tau)}.
	\] 
\end{theorem}

We should note that this simultaneously generalizes both the known Ehrhart series results for order polytopes and $\s$-lecture hall simplices. 
Specifically, if $\s(i)=1$ for all $i$, then we recover the standard order polytope $O(P)$ and the Ehrhart series given by Stanley. 
Given any choice of $\s$, if $P$ is chosen to be  the totally ordered chain, we recover $\Pbn$ and an alternative description of its $h^\ast$-polynomial.

\commentout{
When assuming certain conditions on $(P,\s)$, one can make palindromic statements. 
Suppose that $P$ is a labeled poset, let $\mathcal{E}(P)$ denote the set of covering relations of $P$, and define the sign function $\epsilon:E(P)\to\{-1,1\}$
	\[
	\epsilon(x,y):=\begin{cases}
	1, & \mbox{ if } x<y, \mbox{ and }\\
	-1, & \mbox{ if } x>y.
	\end{cases}
	\]
We say that $P$ is \emph{sign-graded or rank $r$} if
	\[
	\sum_{i=1}^k\epsilon(x_{i-1},x_i)=r
	\] 
for each maximal chain $x_0\prec x_1\prec \cdots \prec x_k$ is $P$.
This definition leads to a well-defined \emph{rank function} on a sign-graded poset, namely
	\[
	\rho(x)=\sum_{i=1}^k\epsilon(x_{i-1},x_i)
	\] 	
where $x_0\prec x_1\prec \cdots \prec x_k=x$ is any saturated chain in $P$ with $x_0$ a minimal element.
A labeled poset is called \emph{sign-ranked} if for each maximal element $x$ of $P$, the subposet $\{y\in P \, : y\preceq x\}$ is sign-graded.
We can now state the following theorem of Br\"and\'{e}n and Leander for such posets.
}

When assuming certain conditions on $(P,\s)$, one can make palindromic statements. 
We say that a poset $P=([n],\preceq)$ is \emph{graded} if it is equipped with a rank function $\rho:[n]\to \Z$ that has the properties that $\rho(x)<\rho(y)$ whenever $x\prec y$ and $\rho(y)=\rho(x)+1$ if $y$ covers  $x$ in $P$.
We can now state the following theorem of Br\"and\'{e}n and Leander for such posets.

\begin{theorem}[\cite{BrandenLeander-PPartition}]\label{signranked}
Suppose that $P=([n],\preceq)$ is a naturally-labeled, graded poset with nonnegative rank function $\rho$ and $\s=\rho+1$.
Then $h^\ast(O(P,\s);z)$ is a palindromic polynomial of degree $n-1$.
Equivalently, $O(P,\s)$ is a Gorenstein polytope of index 2.
\end{theorem}

\subsection{Unimodality and real-rootedness}

Given that the $h^\ast$-polynomials for $O(P,\s)$ can be combinatorially described, it is natural to consider when these polynomials are unimodal. 
Additionally, since all of $O(P,\s)$ are known to have the IDP and since there is a natural family of Gorenstein $O(P,\s)$ arising from sign-ranked posets,  unimodality questions are especially on interest in terms of the Hibi and Ohsugi Conjecture (see Conjecture \ref{conj:HibiOhsugi}) and the question of Schepers and Van Langenhoven (See Question \ref{ques:IDPUnimodal}).

A powerful technique for showing unimodality of polynomials with positive coefficients is showing real-rootedness.
The following two results of Br\"and\'{e}n and Leander show real-rootedness in certain cases. 

\begin{theorem}[\cite{BrandenLeander-PPartition}]
Let $P,Q$ be labeled posets on $[p]$ and $[q]$ respectively and let $\s_P:[p]\to\Zpos$ and $\s_Q:[q]\to \Zpos$.
Let $P\sqcup Q$ denote the disjoint union of $P$ and $Q$ and let $\s_{P\sqcup Q}$ be the unique function on $P\sqcup Q$ which agrees with $\s_P$ and $\s_Q$ on  the components $P$ and $Q$ respectively.
If $h^\ast(O(P,\s_P);z)$ and $h^\ast(O(Q,\s_Q);z)$ are real-rooted, then $h^\ast(O(P\sqcup Q, \s_{P\sqcup Q});z)$ is real-rooted.
\end{theorem}

\begin{theorem}[\cite{BrandenLeander-PPartition}]
Let $P=A_{p_1}\oplus\cdots\oplus A_{p_m}$ be an ordinal sum of anti-chains and let $\s:P\to \Zpos$ be a function which is constant on $A_{p_i}$ for $i\leq i\leq m$. Then $h^\ast(O(P,\s);z)$ is real-rooted.
\end{theorem}

The first theorem follows quickly as it can be shown that  
	\[
	h^\ast(O(P\sqcup Q, \s_{P\sqcup Q};z)=h^\ast(O(P,\s_P);z)\cdot h^\ast(O(Q,\s_Q);z).
	\]
The second theorem requires more sophisticated techniques. In particular, they show that $h^\ast(O(P,\s);z)$ can be expressed as a sum of related polynomials which are an \emph{interlacing sequence}.
We should note that this is similar to the method used by Savage and Visontai \cite{SavageVisontai} to prove Theorem \ref{thm:sEulerRealRoot}.
 The reader should consult  \cite{Branden-Survey} for background and details on this method of proving real-rootedness.

In the special case that $P$ is a naturally-labeled poset with a unique minimal element and nonnegative rank function and where we let $\s=\rho+1$, which is a subcase of that covered by Theorem \ref{signranked},  unimodality can be shown using triangulation results. In particular, the following theorem of Gustafsson and Solus follows with some work from Theorem \ref{ordertriangulation} in the next subsection.

\begin{theorem}[{\cite{GustafssonSolus}}]
Suppose that $P$ is a naturally-labeled, graded poset with a unigue minimal element, a nonnegative rank function $\rho$, and $\s=\rho+1$. Then $h^\ast(O(P,\s);z)$ is unimodal.
\end{theorem}

\subsection{Box unimodal triangulations and local $h^\ast$-polynomials}
To conclude our discussion of results on $\s$-lecture hall order polytopes, we will consider work regarding  box unimodal triangulations and local Ehrhart theory.
This work is due entirely to Gustafsson and Solus \cite{GustafssonSolus}. 
The \emph{$\s$-canonical triangulation} of $O(P,\s)$ is the analogue to the canonical triangulation of the order polytope $O(P)$ introduced by Stanley \cite{Stanley-PosetPolytopes} with the primary difference being that rather than decomposing the polytope into unimodular simplicies indexed by linear extensions of $P$, $O(P,\s)$ is decomposed into $\s$-lecture hall simplices, which are in general not unimodular, in accordance with linear extensions of $P$.
We should note that both the canonical and $\s$-canonical triangulations are regular.
We begin with a theorem on the existence of box unimodal triangulations.

\begin{theorem}[{\cite[Theorem 5.6]{GustafssonSolus}}]\label{ordertriangulation}
Let $O(P,\s)$ be an $\s$-lecture hall order polytope. 
Then the $\s$-canonical triangulation of $O(P,\s)$ is box unimodal. 
Moreover, if $O(P,\s)$ is reflexive, then the $h^\ast$-polynomial $h^\ast(O(P,\s);z)$ is unimodal. 
\end{theorem}

This theorem combined with Theorem \ref{theorem:IDPOrder} provides further support of Conjecture \ref{conj:HibiOhsugi}.
While this statement does not guarantee the existence of a regular, unimodular triangulation, which is the case in the order polytope case, it is nevertheless a very strong statement which gives some evidence for $h^\ast$-polynomial unimodality results in generality.

Given that the local $h^\ast$-polynomials of $\s$-lecture hall simplices are both real-rooted and $\gamma$-positive, it is a natural extension to consider such properties of the local $h^\ast$-polynomials of $O(P,\s)$. 
In general, these strong conditions cannot currently be verified.
However, Theorem \ref{ordertriangulation} can be leveraged, along with results of Katz and Stapledon \cite{KatzStapledon-Ehrhart}, to show the following unimodality result.

\begin{theorem}[{\cite[Theorem 5.9]{GustafssonSolus}}]
The local $h^\ast$-polynomial of an $\s$-lecture hall order polytope is unimodal.
\end{theorem}

\section{Open Problems}\label{sec:Open}
We conclude this survey with a collection of conjectures and questions which have accumulated in the literature, so that the majority of open problems regarding the polyhedral geometry of lecture hall partitions appear in a single location.
We will follow the order of the survey by first addressing open problems for $\s$-lecture hall cones, then $\s$-lecture hall simplices, and culminating with $\s$-lecture hall order polytopes. 

\subsection{Lecture hall cones}

Given that the  Gorenstein property for $\Cc_n^{(\s)}$ has been completely characterized, a natural extension is to consider the level property.
This generates the following question:
	\begin{question}
	Given $\s\in\Zpos^n$, what conditions on $\s$ guarantee that $\Cc_n^{(\s)}$ is level?
	\end{question} 
For this question, there is a notable obstruction; there is no canonical choice of grading for $\Cc_n^{(\s)}$. 
Choices for grading which have been used in the literature include the standard ``height" grading used in the polytope case, as well as the difference of last two coordinates which was used heavily in the case of $\Cc_n^{(1,2,\ldots,n)}$.
After choosing a grading, it seems plausible that one may be able to adapt the strategy of Kohl and the author in the $\Pbn$ case \cite{KohlOlsen-Level}, though one would need to understand the Hilbert basis for any given cone considered. 

Additionally, there are many potential extensions of the Hilbert basis classifications. 
This leads to the following problem and question.
\begin{problem}\label{prob:hilb}
Identify the the Hilbert basis elements of additional combinatorial families of lecture hall cones.
\end{problem} 
\begin{question}\label{ques:hilb}
Is there a reasonable description of the Hilbert bases of $\u$-generated Gorenstein lecture hall cones of dimension 5?
\end{question}
For Problem \ref{prob:hilb}, a possible candidate would be to consider the $(k,\ell)$-sequences, or even a special case such as the $(4,1)$-sequences. The reader should consult \cite{Savage-LHP-Survey} for an overview of these sequences and the deep mathematics behind them.
Regarding Question \ref{ques:hilb}, the Hilbert bases for $\u$-generated Gorenstein lecture hall cones have been classified for dimension $\leq 4$. 
Continuing the classification to higher dimensions should certainly be possible, though there are likely many subcases and enumeration of the cardinality of the Hilbert basis may be significantly more challenging. 

For the special case of $\Cc_n^{(1,2,\ldots, n)}$ where a triangulation of the cone is known, the following conjecture stated by Beck et al. would provide additional geometric information about the cone.

\begin{conjecture}[\cite{BeckEtAl-TriangulationsLHC}]
There exists a regular, flag, unimodular triangulation of $\mathscr{R}_n$ that admits a shelling order
such that the maximal simplices of the triangulation are indexed by $\pi\in\mathfrak{S}_{n-1}$ and each such simplex is
attached along $\des(\pi)$ many of its facets. 
\end{conjecture} 

\subsection{Lecture hall simplices}
While many results are known for the $\s$-lecture hall simplices, there are several open questions in well-studied areas such as triangulations and the levelness property.
To begin, we consider a conjecture stated by Hibi, Tsuchiya, and the author.
\begin{conjecture}[\cite{HibiOlsenTsuchiya-LHPGorensteinIDP}]
Let $\s\in\Zpos^n$ be given. Then $\Pbn$ has a regular, unimodular triangulation.
\end{conjecture}
This has been shown in very particular cases as detailed in the survey. 
However, there is no known counterexample in other cases.
Verifying this conjecture is also crucial for showing results regarding the $\s$-canonical triangulation of   $\s$-lecture hall order polytopes.

While the level property is completely classified in terms of $\s$-inversion sequences, the characterization is often difficult to apply.
Subsequently, it is of interest to determine alternative characterizations for restricted families of $\s$-sequences or to determine conditions which are sufficient though perhaps not necessary for levelness.
This motivates the following question and conjecture.
\begin{question}
Is there an alternative characterization of levelness if $\s$ is a strictly monotone sequence? 
\end{question}
\begin{conjecture}[\cite{KohlOlsen-Level}]
Let $\s\in\Zpos^n$ be a sequence such that there exists some $\c\in\Z^n$ satisfying
	\[
	c_js_{j-1}=c_{j-1}s_j +\gcd(s_{j-1},s_j)
	\]
for $j>1$ with $c_1=1$. Then $\Pbn$ is level.
\end{conjecture}

While there are some results on the Ehrhart positivity for $\s$-lecture hall simplices, the property is largely unstudied.
This motivates the following question.

\begin{question}
What conditions on $\s$ ensure that $\Pbn$ is Ehrhart positive? What conditions on $\s$, besides those in Theorem \ref{thm:EhrhartPos}, ensure that $\Pbn$ is not Ehrhart positive?
\end{question}  

One property which has yet to be studied regarding the $\s$-lecture hall simplices is \emph{reflexive dimension}. 
Given a polytope $\Pc$, the \emph{reflexive dimension} of $\Pc$ is the smallest $m$ such that $\Pc$ is the face of an $m$-dimensional reflexive polytope.
By a result of Haase and Melnikov \cite{HaaseMelnikov}, this is known to be finite for any $\Pc$. This yields the following question.
\begin{question}
Let $\s\in\Zpos^n$ be given. What is the reflexive dimension of $\Pbn$?
\end{question}

\subsection{Lecture hall order polytopes}
We begin the discussion of open problems for $\s$-lecture hall polytopes by mentioning two conjectures of Br\"and\'{e}n and Leander \cite{BrandenLeander-PPartition} in the seminal paper.

\begin{conjecture}[\cite{BrandenLeander-PPartition}]
\label{conj:gammapos}
Suppose $P$ is a naturally-labeled, graded poset  with nonnegative rank function $\rho$ and $\s=\rho+1$, then $h^\ast(O(P,\s);z)$ is $\gamma$-nonnegative. 
\end{conjecture}

\begin{conjecture}[\cite{BrandenLeander-PPartition}]
\label{conj:ordertriangulation}
Suppose that $P$ is a naturally-labeled, graded poset with nonnegative rank function $\rho$ and $\s=\rho+1$, then $O(P,\s)$ has a regular, unimodular triangulation.
\end{conjecture} 

Both of these conjectures are primarily of interest because they imply that $h^\ast(O(P);z)$ is unimodal.
Gustafsson and Solus \cite{GustafssonSolus} proved a slightly weaker unimodality result to that of Conjecture \ref{conj:gammapos} in the case that $P$ has a unique minimal element as discussed previously. 
However, this question is still open in the general setting.
Conjecture \ref{ordertriangulation} is also open in generality. In fact, this is only known to hold in very specific cases, such as the case where $P$  is the totally ordered chain with minimal element $x_0$ and $\rho(x_0)\in\{0,1\}$ where it overlaps with known results for lecture hall simplices.

We conclude with three general questions regarding $\s$-lecture hall order polytopes.

\begin{question}
\label{ques:Gorenstein}
Let $P=([n],\preceq)$ be a labeled poset and $\s:[n]\to \Zpos$. Can we determine conditons of $P$ and $\s$ different than those given in Theorem \ref{signranked} which guarantee that $O(P,\s)$ is Gorenstein?
\end{question}

\begin{question}
\label{ques:level}
Let $P=([n],\preceq)$ be a labeled poset and $\s:[n]\to \Zpos$. What conditions on $P$ and $\s$ guarantee that $O(P,\s)$ is level?
\end{question}

\begin{question}
\label{ques:triangulation}
Let $P=([n],\preceq)$ be a labeled poset and $\s:[n]\to \Zpos$. When does $O(P,\s)$ admit a regular, unimodular triangulation? 
\end{question}

Regarding Questions \ref{ques:Gorenstein} and \ref{ques:level}, characterizations for the Gorenstein and level properties are known for $\s$-lecture hall simplices as discussed in Section \ref{sec:Simplices}, as well as for order polytopes (see \cite{Hibi-LatticesAndStraighteningLaws} for the Gorenstein property and see \cite{HaaseKohlTsuchiya,Miyazaki} for levelness).
It seems reasonable that one may be able to adapt the  methods and techniques to find suitable characterizations in this context. 
One possible approach to Question \ref{ques:triangulation} is to first consider the case of lecture hall simplices and see if the $\s$-canonical triangulation can be exploited.

\bibliographystyle{plain}

\end{document}